\documentclass[12pt]{article}

\usepackage{graphicx,psfrag}
\usepackage{psfrag,eepic,epsfig}
\usepackage{sectsty}
\usepackage{setspace}
\usepackage{color}

\usepackage{setspace}
\usepackage{amsmath, epsfig, cite, ifpdf}
\usepackage{amssymb}
\usepackage{amsfonts}
\usepackage{latexsym}
\usepackage{graphicx,psfrag,eepic,rawfonts}
\usepackage{enumerate}
\usepackage{indentfirst}
\input{epsf}

\newtheorem{thm}{Theorem}[section]
\newtheorem{prop}[thm]{Proposition}

\newtheorem{cor}[thm]{Corollary}

\newtheorem{lem}[thm]{Lemma}

\newcommand{\qed}{{\hfill\rule{4pt}{7pt}}}
\def\pf{\noindent {\it Proof.} }

\numberwithin{equation}{section}

\makeatletter \@addtoreset{equation}{section} \makeatother

\title {\bf Lower bounds of the skew spectral radii and skew energy
of oriented graphs\footnote{Supported by NSFC and PCSIRT. }}

\author{
{\small Xiaolin Chen, Xueliang Li, Huishu Lian}\\
{\small Center for Combinatorics and LPMC-TJKLC}\\
{\small Nankai University, Tianjin 300071, P.R. China}\\
{\small E-mail: chxlnk@163.com; lxl@nankai.edu.cn; lhs6803@126.com}
   }
\date{}
\begin{document}

\maketitle

\begin{abstract}
Let $G$ be a graph with maximum degree $\Delta$, and let
$G^{\sigma}$ be an oriented graph of $G$ with skew adjacency matrix
$S(G^{\sigma})$. The skew spectral radius $\rho_s(G^{\sigma})$ of
$G^\sigma$ is defined as the spectral radius of $S(G^\sigma)$. The
skew spectral radius has been studied, but only few results about
its lower bound are known. This paper determines some lower bounds
of the skew spectral radius, and then studies the oriented graphs
whose skew spectral radii attain the lower bound $\sqrt{\Delta}$.
Moreover, we apply the skew spectral radius to the skew energy of
oriented graphs, which is defined as the sum of the norms of all the
eigenvalues of $S(G^\sigma)$, and denoted by
$\mathcal{E}_s(G^\sigma)$. As results, we obtain some lower bounds
of the skew energy, which improve the known lower bound obtained
by Adiga et al.\\

\noindent\textbf{Keywords:} oriented graph, skew adjacency
matrix, skew spectral radius, skew energy\\

\noindent\textbf{AMS Subject Classification 2010:} 05C20, 05C50,
15A18, 05C90
\end{abstract}

\section{Introduction}

The spectral radius of a graph is one of the fundamental subjects in
spectral graph theory, which stems from the spectral radius of a
matrix. Let $M$ be a square matrix. Then the spectral radius of $M$,
denoted by $\rho(M)$, is defined as the maximum norm of its all
eigenvalues. If $G$ is a simple undirected graph with adjacency
matrix $A(G)$, then the spectral radius of $G$ is defined to be the
spectral radius of $A(G)$, denoted by $\rho(G)$. It is well-known
that $\rho(G)$ is the largest eigenvalue of $A(G)$. The spectral
radius of undirected graphs has been studied extensively and deeply.
For the bounds of the spectral radius, many results have been
obtained. The spectral radius of a graph is related to the chromatic
number, independence number and clique number of the graph; see
\cite{BH} for details. Moreover, the spectral radius of a graph has
applications in graph energy \cite{LSG}.

Recently, the spectral radii of skew adjacency matrices of oriented
graphs have been studied in \cite{CCF,CXZ,SS,Xu,XG}. Let $G^\sigma$
be an oriented graph of $G$ obtained by assigning to each edge of
$G$ a direction such that the induced graph $G^\sigma$ becomes a
directed graph. The graph $G$ is called the underlying graph of
$G^\sigma$. The skew adjacency matrix of $G^\sigma$ is the $n\times
n$ matrix $S(G^\sigma)=(s_{ij})$, where $s_{ij}=1$ and $s_{ji}=-1$
if $\langle v_i,v_j\rangle$ is an arc of $G^\sigma$,
$s_{ij}=s_{ji}=0$ otherwise. It is easy to see that $S(G^\sigma)$ is
a skew symmetric matrix. Thus all the eigenvalues of $S(G^\sigma)$
are pure imaginary numbers or $0's$, which are said to be the skew
spectrum $Sp_s(G^\sigma)$ of $G^\sigma$. The skew spectral radius of
$G^\sigma$, denoted by $\rho_s(G^\sigma)$, is defined as the
spectral radius of $S(G^\sigma)$.

There are only few results about the skew spectral radii of oriented
graphs. Xu and Gong \cite{XG} and Chen et al. \cite{CXZ} studied the
oriented graphs with skew spectral radii no more than $2$. Cavers et
al. \cite{CCF} and Xu \cite{Xu} independently deduced an upper
bound, that is, the skew spectral radius of $G^\sigma$ is dominated
by the spectral radius of the underlying graph $G$. In Section 2, we
further investigate the skew spectral radius of $G^\sigma$, and
obtain some lower bounds for $\rho_s(G^\sigma)$. One of the lower
bounds is based on the maximum degree of the underlying graph $G$,
that is, $\rho_s(G^\sigma)\geq\sqrt{\Delta}$. We, in Section 3,
study the oriented graphs whose skew spectral radii attain the lower
bound $\sqrt{\Delta}$. In Section 4, we apply the skew spectral
radius to the skew energy of oriented graphs.

The skew energy $\mathcal{E}_s(G^\sigma)$ of $G^\sigma$, first
introduced by Adiga et al. \cite{ABC}, is defined as the sum of the
norms of all the eigenvalues of $S(G^\sigma)$. They obtained that
for any oriented graph $G^\sigma$ with $n$ vertices, $m$ arcs and
maximum degree $\Delta$,
\begin{equation}\label{Skew}
\sqrt{2m+n(n-1)\left(\det(S)\right)^{2/n}}\leq\mathcal{E}_s(G^\sigma)\leq n\sqrt{\Delta},
\end{equation}
where $S$ is the skew adjacency matrix of $G^\sigma$ and $\det(S)$
is the determinant of $S$.

The upper bound that $\mathcal{E}_s(G^\sigma)=n\sqrt{\Delta}$ is
called the optimum skew energy. They further proved that an oriented
graph $G^\sigma$ has the optimum skew energy if and only if its skew
adjacency matrix $S(G^\sigma)$ satisfies that
$S(G^\sigma)^TS(G^\sigma)=\Delta I_n$, or equivalently, all
eigenvalues of $G^\sigma$ are equal to $\sqrt{\Delta}$, which
implies that $G$ is a $\Delta$-regular graph.

From the discussion of Section $3$, it is interesting to find that
if the underlying graph is regular, the oriented graphs with
$\rho_s(G^\sigma)=\sqrt{\Delta}$ have the optimal skew energy.
Moreover, by applying the lower bounds of the skew spectral radius
obtained in Section $2$, we derive some lower bounds of the skew
energy of oriented graphs, which improve the lower bound in
(\ref{Skew}).

Throughout this paper, when we simply mention a graph, it means a
simple undirected graph. When we say the maximum degree, average
degree, neighborhood, etc. of an oriented graph, we mean the same as
those in its underlying graph, unless otherwise stated. For
terminology and notation not defined here, we refer to the book of
Bondy and Murty \cite{BM}.

\section{Lower bounds of the skew spectral radius}

In this section, we deduce some lower bounds of the skew spectral
radii of oriented graphs, which implies the relationships between
skew spectral radius and some graph parameters.

We begin with some definitions. Let $G^\sigma$ be an oriented graph
of a graph $G$ with vertex set $V$. Denote by $S(G^\sigma)$ and
$\rho_s(G^\sigma)$ the skew adjacency matrix and the skew spectral
radius of $G^\sigma$, respectively. For any two disjoint subsets
$A,B$ of $V$, we denote by $\gamma(A,B)$ the number of arcs whose
tails are in $A$ and heads in $B$. Let $G^\sigma[A]$ be the subgraph
of $G^\sigma$ induced by $A$, where $G^\sigma[A]$ has vertex set $A$
and contains all arcs of $G^\sigma$ which join vertices of $A$. Then
we can deduce a lower bound of $\rho_s(G^\sigma)$ as follows.
\begin{thm}\label{lower}
Let $G$ be a graph with vertex set $V$, and let $G^\sigma$ be an
oriented graph of $G$ with skew adjacency matrix $S(G^\sigma)$ and
skew spectral radius $\rho_s(G^\sigma)$. Then for any two nonempty subsets
$A,B\subseteq V$,
\begin{equation*}
\rho_s(G^\sigma)\geq\frac{|\gamma(A,B)-\gamma(B,A)|}{\sqrt{|A||B|\,}}\,,
\end{equation*}
where $|A|$ and $|B|$ are the number of elements of $A$ and $B$, respectively.
\end{thm}
\pf Suppose that $|A|=k$ and $|B|=l$, and let $C=A\cap B$ with order $t$, $A_1=A-C$ with order $k-t$,
$B_1=B-C$ with order $l-t$, and $D=V-A\cup B$ with order $n+t-k-l$.
With suitable labeling of the vertices of $V$, the skew adjacency matrix
$S$ can be formulated as follows:
$$\left(\begin{array}{cccc} S_{11} & S_{12} & S_{13} & S_{14}\\ -S_{12}^T & S_{22} & S_{23} & S_{24}\\
-S_{13}^T & -S_{23}^T & S_{33} & S_{34}\\ -S_{14}^T & -S_{24}^T & -S_{34}^T & S_{44} \end{array}\right),$$
where $S_{11}=S(G^\sigma[A_1])$ is the skew adjacency matrix of the induced
oriented graph $G^\sigma[A_1]$ with order $k-t$, $S_{22}=S(G^\sigma[C])$
is the skew adjacency matrix of $G^\sigma[C]$ with order $t$ and
$S_{33}=S(G^\sigma[B_1])$ is the skew adjacency matrix of
$G^\sigma[B_1]$ with order $l-t$, $S_{44}=S(G^\sigma[D])$ is the skew adjacency matrix of the induced
oriented graph $G^\sigma[D]$ with order $n+t-k-l$.

Let $H=(-i)S$. Since $S$ is skew symmetric, $H$ is an Hermitian
matrix. By Rayleigh-Ritz theorem, we obtain that
\begin{equation*}
\rho_s(G^\sigma)=\rho(S)=\rho(H)=\max_{x\in \mathbb{C}^n}\frac{x^*Hx}{x^*x},
\end{equation*}
where $\mathbb{C}^n$ is the complex vector space of dimension $n$
and $x^*$ is the conjugate transpose of $x$. We take
$x=\Big(\underbrace{\frac{1}{\sqrt{k}},\ldots,\frac{1}{\sqrt{k}}}_{k-t},
\underbrace{\frac{1}{\sqrt{k}}+\frac{i}{\sqrt{l}},\ldots,\frac{1}{\sqrt{k}}+\frac{i}{\sqrt{l}}}_{t},
\underbrace{\frac{i}{\sqrt{l}},\ldots,\frac{i}{\sqrt{l}}}_{l-t},
0,\ldots,0\Big)^T=\left(\frac{1}{\sqrt{k}}\mathbf{1}_{k-t}^T,(\frac{1}{\sqrt{k}}+\frac{i}{\sqrt{l}})\mathbf{1}_t^T,
\frac{i}{\sqrt{l}}\mathbf{1}_{l-t}^T,\mathbf{0}^T\right)^T$ in the above equality and derive that
\begin{eqnarray*}
\rho_s(G^\sigma)&\geq &\frac{x^*Hx}{x^*x}=\frac{-i}{2}\left(\frac{1}{k}\mathbf{1}_{k-t}^TS_{11}\mathbf{1}_{k-t}
+\left(\frac{1}{k}+\frac{1}{l}\right)\mathbf{1}_t^TS_{22}\mathbf{1}_t+\frac{1}{l}\mathbf{1}_{l-t}^TS_{33}\mathbf{1}_{l-t}\right.\\
&+&\left(\frac{1}{k}+\frac{i}{\sqrt{kl}}\right)\mathbf{1}_{k-t}^TS_{12}\mathbf{1}_{t}-
\left(\frac{1}{k}-\frac{i}{\sqrt{kl}}\right)\mathbf{1}_{t}^TS_{12}^T\mathbf{1}_{k-t}\\
&+&\frac{i}{\sqrt{kl}}\mathbf{1}_{k-t}^TS_{13}\mathbf{1}_{l-t}-\frac{-i}{\sqrt{kl}}\mathbf{1}_{l-t}^TS_{13}^T\mathbf{1}_{k-t}\\
&+&\left.\left(\frac{i}{\sqrt{kl}}+\frac{1}{l}\right)\mathbf{1}_{t}^TS_{23}\mathbf{1}_{l-t}-
\left(\frac{-i}{\sqrt{kl}}+\frac{1}{l}\right)\mathbf{1}_{l-t}^TS_{23}^T\mathbf{1}_{t}\right).
\end{eqnarray*}

Note that
$\mathbf{1}_{k-t}^TS_{11}\mathbf{1}_{k-t}=\mathbf{1}_t^TS_{22}\mathbf{1}_t=\mathbf{1}_{l-t}^TS_{33}\mathbf{1}_{l-t}=0$,
since $S_{11}$, $S_{22}$ and $S_{33}$ are both skew symmetric. Moreover, it
can be verified that
\begin{equation*}
\mathbf{1}_{k-t}^TS_{12}\mathbf{1}_{t}=\mathbf{1}_{t}^TS_{12}^T\mathbf{1}_{k-t}=\gamma(A_1,C)-\gamma(C,A_1),
\end{equation*}
\begin{equation*}
\mathbf{1}_{k-t}^TS_{13}\mathbf{1}_{l-t}=\mathbf{1}_{l-t}^TS_{13}^T\mathbf{1}_{k-t}=\gamma(A_1,B_1)-\gamma(B_1,A_1),
\end{equation*}
and
\begin{equation*}
\mathbf{1}_{t}^TS_{23}\mathbf{1}_{l-t}=\mathbf{1}_{l-t}^TS_{23}^T\mathbf{1}_{t}=\gamma(C,B_1)-\gamma(B_1,C).
\end{equation*}
Note that $\gamma(A,B)-\gamma(B,A)=\gamma(A_1,C)+\gamma(A_1,B_1)+\gamma(C,B_1)-\gamma(C,A_1)-\gamma(B_1,A_1)-\gamma(B_1,C)$.

Combining the above equalities, it follows that
\begin{equation*}
\rho_s(G^\sigma)\geq\frac{\gamma(A,B)-\gamma(B,A)}{\sqrt{kl}}.
\end{equation*}

Similarly, we can derive that
$\rho_s(G^\sigma)\geq\frac{\gamma(B,A)-\gamma(A,B)}{\sqrt{kl}}$. The
proof is thus complete.\qed

The above theorem implies a lower bound of $\rho_s(G^\sigma)$ in
terms of the maximum degree of $G$. Before proceeding, it is
necessary to introduce the notion of switching-equivalent \cite{LL}
of two oriented graphs. Let $G^\sigma$ be an oriented graph of $G$
and $W$ be a subset of its vertex set. Denote
$\overline{W}=V(G^\sigma)\setminus W$. Another oriented graph
$G^\tau$ of $G$, obtained from $G^\sigma$ by reversing the
orientations of all arcs between $W$ and $\overline{W}$, is said to
be obtained from $G^\sigma$ by switching with respect to $W$. Two
oriented graphs $G^\sigma$ and $G^\tau$ are called
switching-equivalent if $G^\tau$ can be obtained from $G^\sigma$ by
a sequence of switchings. The following lemma shows that the
switching operation keeps skew spectrum unchanged.
\begin{lem}\cite{LL}\label{switch}
Let $G^\sigma$ and $G^\tau$ be two oriented graphs of a graph $G$.
If $G^\sigma$ and $G^\tau$ are switching-equivalent, then $G^\sigma$
and $G^\tau$ have the same skew spectra.
\end{lem}
\begin{cor}\label{maximum}
Let $G^\sigma$ be an oriented graph of $G$ with maximum degree
$\Delta$. Then
\begin{equation*}
\rho_s(G^\sigma)\geq \sqrt{\Delta}.
\end{equation*}
\end{cor}
\pf Let $G^\tau$ be an oriented graph of $G$ obtained from
$G^\sigma$ by switching with respect to every neighbor of $v$ if
necessary, such that all arcs incident with $v$ have the common tail
$v$. Then $G^\sigma$ and $G^\tau$ are switching-equivalent. By Lemma
\ref{switch}, $G^\sigma$ and $G^\tau$ have the same skew spectra.
Consider the oriented graph $G^\tau$ and let $A=\{v\}$ and $B=N(v)$.
Obviously, $\gamma(A,B)-\gamma(B,A)=\Delta$. By Theorem \ref{lower},
$\rho_s(G^\sigma)=\rho_s(G^\tau)\geq
\frac{\Delta}{\sqrt{\Delta}}=\sqrt{\Delta}$.\qed

It is known \cite{BH} that for any undirected tree $T$, $\rho(T)\geq
\bar{d}$, where $\bar{d}$ is the average degree of $T$. Moreover,
all oriented trees of $T$ have the same skew spectra which are equal
to $i$ times the spectrum of $T$; see \cite{SS}. Therefore, for any
oriented tree $T^\sigma$ of $T$, $\rho_s(T^\sigma)=\rho(T)\geq
\bar{d}$. We next consider a general graph. The following lemma
\cite{AS} is necessary.
\begin{lem}\cite{AS}\label{bipartite}
Let $G=(V,E)$ be a graph with $n$ vertices and $m$ edges. Then $G$
contains a bipartite subgraph with at least $m/2$ edges.
\end{lem}

Then we obtain the following result for a general graph by applying
Theorem \ref{lower}.
\begin{cor}
For any simple graph $G$ with average degree $\bar{d}$, there exists
an oriented graph $G^\sigma$ of $G$ such that
\begin{equation*}
\rho_s(G^\sigma)\geq \frac{\,\bar{d}\,}{2}.
\end{equation*}
\end{cor}
\pf By Lemma \ref{bipartite}, $G$ contains a bipartite subgraph
$H=(A,B)$ with at least $m/2$ edges. We give an orientation of $G$
such that all arcs between $A$ and $B$ go from $A$ to $B$ and the
directions of the other arcs are arbitrary. By Theorem \ref{lower}
and Lemma \ref{bipartite},
\begin{equation*}
\rho_s(G^\sigma)\geq\frac{|\gamma(A,B)-\gamma(B,A)|}{\sqrt{|A||B|}}
\geq\frac{2}{n}\,|\gamma(B,A)-\gamma(A,B)|
\geq\frac{m}{n}=\frac{\bar{d}}{2}.
\end{equation*}\qed

Hofmeister \cite{H} and Yu et al. \cite{YLT} deduced a lower bound
of the spectral radius of a graph in terms of its degree sequence.
Specifically, let $G$ be a connected graph with degree sequence
$d_1,d_2,\ldots,d_n$. Then $\rho(G)\geq
\sqrt{\frac{1}{n}\sum_{i=1}^nd_i^2}$. Similarly, for an oriented
graph, we consider the relation between its skew spectral radius and
vertex degrees, where it should be taken into account the out-degree
and in-degree of every vertex.

Let $G^\sigma$ be an oriented graph with vertex set
$\{v_1,v_2,\ldots,v_n\}$. Denote by $d_i^+$ and $d_i^-$ the
out-degree and in-degree of $v_i$ in $G^\sigma$, respectively. Let
$\tilde{d_i}=d_i^+-d_i^-$. Then we establish a lower bound of the
skew spectral radius of $G^\sigma$ as follows.
\begin{thm}\label{degree}
Let $G^\sigma$ be an oriented graph with vertex set
$\{v_1,v_2,\ldots,v_n\}$ and skew spectral radius
$\rho_s(G^\sigma)$. Then
\begin{equation*}
\rho_s(G^\sigma)\geq\sqrt{\frac{\tilde{d_1}^2+\tilde{d_2}^2+\cdots+\tilde{d_n}^2}{n}}\,.
\end{equation*}
\end{thm}
\pf If $G^\sigma$ is an Eulerian digraph, then the right-hand of the
above inequality is $0$, and the inequality is obviously true since
$\rho_s(G^\sigma)\geq 0$ holds always. So, we always assume that
$G^\sigma$ is not Eulerian in the following. Let $S=[s_{ij}]$ be the
skew adjacency matrix of $G^\sigma$. Then
$\rho_s(G^\sigma)=\rho(S)=\sqrt{\rho(S^TS)}$. We consider the
spectral radius of $S^TS$. Since $S^TS$ is symmetric, by
Rayleigh-Ritz theorem,
\begin{equation*}
\rho(S^TS)=\max_{x\in \mathbb{R}^n}\frac{x^T(S^TS)x}{x^Tx}.
\end{equation*}
We take $x=\frac{1}{\sqrt{\tilde{d_1}^2+\tilde{d_2}^2+\cdots+\tilde{d_n}^2}}\left(\tilde{d_1},
\tilde{d_2},\ldots,\tilde{d_n}\right)^T$ in the above equality and obtain that
\begin{equation*}
\rho(S^TS)\geq x^T(S^TS)x=(Sx)^T(Sx).
\end{equation*}
It is easy to compute that
\begin{equation*}
Sx=\frac{1}{\sqrt{\tilde{d_1}^2+\tilde{d_2}^2+\cdots+\tilde{d_n}^2}}\left(\sum_{j=1}^ns_{1j}\tilde{d_j},
\sum_{j=1}^ns_{2j}\tilde{d_j},\ldots,\sum_{j=1}^ns_{nj}\tilde{d_j}\right)^T.
\end{equation*}
Applying the Cauchy-Schwarz's inequality, we obtain that
\begin{eqnarray*}
(Sx)^TSx&=&\frac{1}{\tilde{d_1}^2+\tilde{d_2}^2+\cdots+\tilde{d_n}^2}
\left[\left(\sum_{j=1}^ns_{1j}\tilde{d_j}\right)^2+
\left(\sum_{j=1}^ns_{2j}\tilde{d_j}\right)^2+\cdots+
\left(\sum_{j=1}^ns_{nj}\tilde{d_j}\right)^2\right]\\
&\geq& \frac{n}{\tilde{d_1}^2+\tilde{d_2}^2+\cdots+\tilde{d_n}^2}
\left(\frac{\sum_{j=1}^ns_{1j}\tilde{d_j}+
\sum_{j=1}^ns_{2j}\tilde{d_j}+\cdots+\sum_{j=1}^ns_{nj}\tilde{d_j}}{n}\right)^2.
\end{eqnarray*}
Note that
\begin{equation*}
\sum_{j=1}^ns_{1j}\tilde{d_j}+
\sum_{j=1}^ns_{2j}\tilde{d_j}+\cdots+\sum_{j=1}^ns_{nj}\tilde{d_j}=(1,1,\ldots,1)Sx=-\tilde{d_1}^2
-\tilde{d_2}^2-\cdots-\tilde{d_n}^2.
\end{equation*}
Therefore,
\begin{equation*}
(Sx)^T(Sx)\geq \frac{\tilde{d_1}^2+\tilde{d_2}^2+\cdots+\tilde{d_n}^2}{n}.
\end{equation*}
We thus conclude that
\begin{equation*}
\rho_s(G^\sigma)\geq\sqrt{(Sx)^T(Sx)}\geq \sqrt{\frac{\tilde{d_1}^2+\tilde{d_2}^2+\cdots+\tilde{d_n}^2}{n}}.
\end{equation*}
This completes the proof.\qed

\noindent\textbf{Remark 2.1} Theorem \ref{degree} can also implies
Corollary \ref{maximum}.

\noindent {\bf \emph{Another proof of Corollary \ref{maximum}}}.
Suppose that $d_G(v_1)=\Delta$ and
$N(v_1)=\{v_2,v_3,\ldots,v_{\Delta+1}\}$. Let $G^\tau$ be an
oriented graph of $G$ obtained from $G^\sigma$ by a sequence of
switchings such that all arcs incident with $v_1$ have the common
tail $v_1$. Then by Lemma \ref{switch},
$\rho_s(G^\sigma)=\rho_s(G^\tau)$. Let $H^\tau$ be the subgraph of
$G^\tau$ induced by $v_1$ and its all adjacent vertices. Let
$H_1=(-i)S(G^\tau)$ and $H_2=(-i)S(H^\tau)$. Then
$\rho_s(G^\tau)=\rho(H_1)$ and $\rho_s(H^\tau)=\rho(H_2)$. Note that
$H_1$ and $H_2$ are both Hermitian matrices. By interlacing of
eigenvalues (\cite{HJ}, 4.3.16 Corollary),
$\rho(H_1)=\lambda_1(H_1)\geq\lambda_1(H_2)=\rho(H_2)$.

Suppose that for any $1\leq i\leq \Delta$, the vertex $v_i$ has
out-degree $t_i^+$ and in-degree $t_i^-$ in $H^\tau$. Let
$\tilde{t_i}=t_i^+-t_i^-$. It can be found that
$\tilde{t_1}+\tilde{t_2}+\cdots+\tilde{t}_{\Delta+1}=0$ and
$\tilde{t_1}=\Delta$. It follows that
$\tilde{t_1}^2+\tilde{t_2}^2+\cdots+\tilde{t}_{\Delta+1}^2\geq
\Delta^2+\Delta$.  Then by Theorem \ref{degree},
\begin{equation*}
\rho_s(H^\tau)\geq\sqrt{\frac{\tilde{t_1}^2+\tilde{t_2}^2+\cdots+\tilde{t}_{\Delta+1}^2}{\Delta+1}}
\geq\sqrt{\Delta}.
\end{equation*}
Now we conclude that
$\rho_s(G^\sigma)=\rho_s(G^\tau)\geq\rho_s(H^\tau)\geq\sqrt{\Delta}$,
which implies Corollary \ref{maximum}.\qed

It is known that Wilf \cite{Wilf} considered the relation between
the spectral radius and the chromatic number of a graph. As for the
oriented graphs, Sopena \cite{Sopena} introduced the notion of
oriented chromatic number, which motivates us to consider the
relation between the skew spectral radius and the oriented chromatic
number of an oriented graph.

Let $G^\sigma$ be an oriented graph with vertex set $V$. An oriented
$k$-coloring of $G^\sigma$ is a partition of $V$ into $k$ color
classes such that no two adjacent vertices belong to the same color
class, and all the arcs between two color classes have the same
direction. The oriented chromatic number of $G^\sigma$, denoted by
$\chi_o(G^\sigma)$, is defined as the smallest number $k$ satisfying
that $G^\sigma$ admits an oriented $k$-coloring. The following
theorem presents a lower bound of the skew spectral radius in terms
of the average degree and oriented chromatic number.
\begin{thm}
Let $G$ be a graph with average degree $\bar{d}$. Let $G^\sigma$ be
an oriented graph of $G$ with skew spectral radius
$\rho_s(G^\sigma)$ and oriented chromatic number $\chi_o(G^\sigma)$.
Then
\begin{equation*}
\rho_s(G^\sigma)\geq \frac{\bar{d}}{\chi_o(G^\sigma)-1}\,.
\end{equation*}
\end{thm}
\pf Suppose that $G^\sigma$ contains $m$ arcs and
$\chi_o(G^\sigma)=k$. Let $\{V_1,V_2,\ldots,V_k\}$ be an oriented
$k$-coloring of $G^\sigma$. Denote
$a_{ij}=|\gamma(V_i,V_j)-\gamma(V_j,V_i)|$. By the definition of an
oriented $k$-coloring, we get that either $\gamma(V_i,V_j)=0$ or
$\gamma(V_j,V_i)=0$. It follows that $\sum_{i<j}a_{ij}=m$. By
Theorem \ref{lower}, we derive that
\begin{eqnarray*}
\rho_s(G^\sigma)&\geq&\max_{i<j}\frac{a_{ij}}{\sqrt{|V_i||V_j|}}
\geq\max_{i<j}\frac{2a_{ij}}{|V_i|+|V_j|}\\
&\geq&\frac{2\sum_{i<j}a_{ij}}{\sum_{i<j}(|V_i|+|V_j|)}
=\frac{2m}{(k-1)n}=\frac{\bar{d}}{\chi_o(G^\sigma)-1}\,.
\end{eqnarray*}
The proof is now complete. \qed

\section{Oriented graphs with skew spectral radius $\sqrt{\Delta}$}

From the previous section, we know that for any oriented graph
$G^\sigma$, $\rho_s(G^\sigma)\geq\sqrt{\Delta}$. In this section, we
investigate the oriented graphs with
$\rho_s(G^\sigma)=\sqrt{\Delta}$.

We first recall the following proposition on the skew spectra of
oriented graphs.
\begin{prop}\label{eigen2}
Let $\{i\lambda_1,i\lambda_2,\ldots,i\lambda_n\}$ be the skew
spectrum of $G^\sigma$, where
$\lambda_1\geq\lambda_2\geq\cdots\geq\lambda_n$. Then (1)
$\lambda_j=-\lambda_{n+1-j}$ for all $1\leq j\leq n$; (2) when $n$
is odd, $\lambda_{(n+1)/2}=0$ and when $n$ is even,
$\lambda_{n/2}\geq 0$; and (3) $\sum_{j=1}^{n}\lambda_j^2=2m$.
\end{prop}

An oriented regular graph is an oriented graph of a regular graph.
We consider the case of oriented regular graphs, which is associated
with optimum skew energy oriented graphs.

\begin{thm}\label{regular}
Let $G^\sigma$ be an oriented graph of a $\Delta$-regular graph $G$
with skew adjacency matrix $S$. Then
$\rho_s(G^\sigma)=\sqrt{\Delta}$ if and only if $S^TS=\Delta I_n$,
i.e., $G^\sigma$ has the optimum skew energy.
\end{thm}
\pf Let $\{i\lambda_1,i\lambda_2,\ldots,i\lambda_n\}$ be the skew
spectrum of $G^\sigma$ with
$\lambda_1\geq\lambda_2\geq\cdots\geq\lambda_n$. By Proposition
\ref{regular}, we get that
$\lambda_1^2+\lambda_2^2+\cdots+\lambda_n^2=2m=n\Delta$ and
$\lambda_1\geq|\lambda_i|$ for any $2\leq i\leq n$. It follows that
$\lambda_1=\rho_s(G^\sigma)$ and
$\lambda_1^2+\lambda_2^2+\cdots+\lambda_n^2\leq
n(\rho_s(G^\sigma))^2$. If $\rho_s(G^\sigma)=\sqrt{\Delta}$, then we
can conclude that
$\lambda_1=|\lambda_2|=\cdots=|\lambda_n|=\sqrt{\Delta}$, that is to
say, $S^TS=\Delta I_n$. The converse implication follows easily.
\qed

The above theorem gives a good characterization for an oriented
$\Delta$-regular graph with $\rho_s(G^\sigma)=\sqrt{\Delta}$, which
says that any two rows and any two columns of its skew adjacency
matrix are all orthogonal. For any oriented graph $G^\sigma$ with
$\rho_s(G^\sigma)=\sqrt{\Delta}$, we also consider the orthogonality
of its skew adjacency matrix and obtain an extended result as
follows.
\begin{thm}\label{general}
Let $G$ be a graph with vertex set $V=\{v_1,v_2,\ldots,v_n\}$ and
maximum degree $\Delta$. Let $G^\sigma$ be an oriented graph of $G$
with $\rho_s(G^\sigma)=\sqrt{\Delta}$ and skew adjacency matrix
$S=(S_1,S_2,\ldots,S_n)$. If $d_G(v_i)=\Delta$, then
$(S_i,S_j)=S_i^TS_j=0$ for any $j\neq i$.
\end{thm}
\pf Without loss of generality, suppose that $d_G(v_1)=\Delta$. It
is sufficient to consider the matrix $S^TS$ and prove that for any
$j\neq 1$, $(S^TS)_{1j}=0$. Notice from
$\rho_s(G^\sigma)=\sqrt{\Delta}$ that $\Delta$ is the maximum
eigenvalue of $S^TS$. Suppose that $\Delta$ is an eigenvalue of
$S^TS$ with multiplicity $l$. Since $S^TS$ is a real symmetric
matrix, there exists an orthogonal matrix $P$ such that
$S^TS=P^TDP$, where $D$ is the diagonal matrix with the form
$D=\text{diag}(\Delta,\ldots,\Delta,u_{l+1},\ldots,u_n)$ with $0\leq
u_i<\Delta$. Denote $P=(P_1,P_2,\ldots,P_n)=(p_{ij})$.

Note that $(S^TS)_{11}=\Delta$ since $d(v_1)=\Delta$. It follows
that $P_1^TDP_1=\Delta$, that is,
\begin{equation*}
\Delta p_{11}^2+\cdots+\Delta p_{l1}^2+u_{l+1}p_{l+1,1}^2+\cdots+u_n p_{n1}^2=\Delta.
\end{equation*}
Since $P_1^TP_1=1$, we derive that
$p_{l+1,1}=p_{l+2,1}=\cdots=p_{n1}=0$. Then for any $j\neq1$, we
compute the $(1,j)$-entry of $S^TS$ as follows.
\begin{eqnarray*}
(S^TS)_{1j}&=&P_1^TDP_j\\
&=&\Delta p_{11}p_{1j}+\cdots+\Delta p_{l1}p_{j1}+u_{l+1}p_{l+1,1}p_{l+1,j}+\cdots+u_n p_{n1}p_{nj}\\
&=&\Delta(p_{11}p_{1j}+\cdots+p_{l1}p_{lj})=\Delta P_1^TP_j=0.
\end{eqnarray*}
The last equality holds due to the orthogonality of $P$. The proof
is now complete.\qed

Comparing Theorem \ref{regular} with Theorem \ref{general}, it is
natural to ask whether the converse of Theorem \ref{general} holds,
that is, whether the condition that $d_G(v_i)=\Delta$ and
$(S_i,S_j)=0$ for every vertex $v_i$ with maximum degree and $j\neq
i$ implies that $\rho_s(G^\sigma)=\sqrt{\Delta}$. In what follows,
we show that it is not always true by constructing a counterexample,
but we can still obtain that $i\sqrt{\Delta}$ is an eigenvalue of
$G^\sigma$.
\begin{thm}\label{eigen}
Let $G$ be a graph with vertex set $V=\{v_1,v_2,\ldots,v_n\}$ and
maximum degree $\Delta$. Let $G^\sigma$ be an oriented graph of $G$
with skew adjacency matrix $S=(S_1,S_2,\ldots,S_n)$. If there exists
a vertex $v_i$ with maximum degree $\Delta$ such that for any $j\neq
i$, $(S_i,S_j)=0$, then $i\sqrt{\Delta}$ is an eigenvalue of
$G^\sigma$.
\end{thm}
\pf Without loss of generality, suppose that $d_G(v_1)=\Delta$.
Since $(S_1,S_j)=0$ for any $j\neq 1$, we obtain that
$SS^T=\left(\begin{array}{cc}\Delta & \mathbf{0}\\ \mathbf{0}^T &
\mathbf{*}\end{array}\right)$. Then $\Delta$ is an eigenvalue of
$SS^T$, which follows that $i\sqrt{\Delta}$ is an eigenvalue of
$G^\sigma$. The proof is thus complete.\qed

\begin{figure}[h,t,b,p]
\begin{center}
\scalebox{1}[1]{\includegraphics{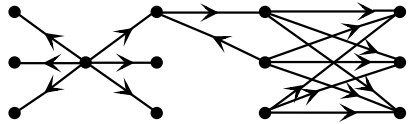}}
\end{center}
\caption{The counterexample $G^\sigma$}\label{Fig1}
\end{figure}

\noindent\textbf{Example 3.1} Let $G^\sigma$ be the oriented graph
depicted in Figure \ref{Fig1}, which has the maximum degree $6$. It
can be verified that $G^\sigma$ satisfies the conditions of Theorem
\ref{eigen} and $\sqrt{6}i$ is an eigenvalue of $G^\sigma$. But we
can compute that $\rho_s(G^\sigma)\approx 3.1260 > \sqrt{6}$.

\section{Lower bounds of the skew energy of $G^\sigma$}

Similar to the McClelland's lower bound for the energy of undirected
graphs, Adiga et al. in \cite{ABC} got a lower bound for the skew
energy of oriented graphs, that is,
$\mathcal{E}_s(G^\sigma)\geq\sqrt{2m+n(n-1)\left(\det(S)\right)^{2/n}}$,
where $S$ is the skew adjacency matrix of $G^\sigma$. This bound is
also called the McClelland's lower bound of skew energy. In this
section, we obtain some new lower bounds for the skew energy of
oriented graphs.

In view of Proposition \ref{eigen2}, we reconsider the McClelland's
lower bound and establish a new lower bound of
$\mathcal{E}_s(G^\sigma)$.
\begin{thm}\label{Mc1}
Let $G^\sigma$ be an oriented graph with $n$ vertices, $m$ arcs and
skew adjacency matrix $S$. Then
\begin{equation}\label{Mc}
\mathcal{E}_s(G^\sigma)\geq\sqrt{4m+n(n-2)\left(\det(S)\right)^{2/n}}.
\end{equation}
\end{thm}
\pf By Proposition \ref{eigen2}, we have
\begin{equation*}
\left(\mathcal{E}_s(G^\sigma)\right)^2=\left(2\sum_{j=1}^{\lfloor n/2\rfloor}|\lambda_j|\right)^2
=4\sum_{j=1}^{\lfloor n/2\rfloor}\lambda_j^2+4\sum_{1\leq i\neq j\leq \lfloor n/2 \rfloor}|\lambda_i||\lambda_j|.
\end{equation*}
If $n$ is odd, $\det(S)=0$ and
$\left(\mathcal{E}_s(G^\sigma)\right)^2 \geq4\sum_{j=1}^{\lfloor
n/2\rfloor}\lambda_j^2=4m$. If $n$ is even, by the
arithmetic-geometric mean inequality, we have that
\begin{equation*}
\left(\mathcal{E}_s(G^\sigma)\right)^2
=4\sum_{j=1}^{\lfloor n/2\rfloor}\lambda_j^2+4\sum_{1\leq i\neq j\leq \lfloor n/2 \rfloor}|\lambda_i||\lambda_j|\geq 4m+n(n-2)(\det(S))^{2/n}.
\end{equation*}
The proof is thus complete.\qed

\noindent\textbf{Remark 4.1} The lower bound in the above theorem is
better than the McClelland's bound. In fact, we find that
\begin{equation}\label{inequ1}
(\det(S))^\frac{1}{n}=\left(\prod_{j=1}^{n}|\lambda_j|\right)^\frac{1}{n}\leq \sqrt{\frac{\sum_{j=1}^{n}\lambda_j^2}{n}}= \sqrt{\frac{2m}{n}}.
\end{equation}
Therefore, we deduce that
\begin{equation*}
4m+n(n-2)\left(\det(S)\right)^\frac{2}{n}\geq 2m+n(n-1)\left(\det(S)\right)^\frac{2}{n}.
\end{equation*}

An oriented graph $G^\sigma$ is said to be singular if $\det(S)=0$
and nonsingular otherwise. In what follows, we only consider the
oriented graphs with $\det(S)\neq 0$. Note that if $G^\sigma$ is
nonsingular, then $n$ must be even and $\det(S)$ is positive. We
next derive a lower bound of the skew energy for nonsingular
oriented graphs in terms of the order $n$, the maximum degree
$\Delta$ and $\det(S)$.

\begin{thm}\label{ELB}
Let $G^\sigma$ be a nonsingular oriented graph with order $n$,
maximum degree $\Delta$ and skew adjacency matrix $S$. Then
\begin{equation}\label{ELB1}
\mathcal{E}_s(G^\sigma)\geq 2\sqrt{\Delta}+(n-2)\left(\frac{\det(S)}{\Delta}\right)^{\frac{1}{n-2}}
\end{equation}
and equality holds if and only if $\lambda_1=\sqrt{\Delta}$ and
$\lambda_2=\cdots=\lambda_{n/2}$.
\end{thm}
\pf Using the arithmetic-geometric mean inequality, we obtain that
\begin{eqnarray*}
\mathcal{E}_s(G^\sigma)&=&\sum_{j=1}^{n}|\lambda_j|
=2\lambda_1+\sum_{j=2}^{n-1}|\lambda_j|\geq 2\lambda_1+(n-2)\left(\prod_{j=2}^{n-2}|\lambda_j|\right)^{\frac{1}{n-2}}\\
&=&2\lambda_1+(n-2)\left(\frac{\det(S)}{\lambda_1^2}\right)^{\frac{1}{n-2}}
\end{eqnarray*}
with equality if and if $|\lambda_2|=\cdots=|\lambda_{n-1}|$.

Let $f(x)=2x+(n-2)\left(\frac{\det(S)}{x^2}\right)^\frac{1}{n-2}$.
Then
$f'(x)=2-2\left(\det(S)\right)^\frac{1}{n-2}x^{-\frac{n}{n-2}}$.

It is easy to see that the function $f(x)$ is increasing for $x\geq \left(\det(S)\right)^\frac{1}{n}$.
By Inequality (\ref{inequ1}), we have $\left(\det(S)\right)^\frac{1}{n}\leq\sqrt{2m/n}\leq \sqrt{\Delta}$.

Therefore,
\begin{equation*}
\mathcal{E}_s(G^\sigma)\geq f(\sqrt{\Delta})=2\sqrt{\Delta}+(n-2)\left(\frac{\det(S)}{\Delta}\right)^\frac{1}{n-2}.
\end{equation*}
and equality holds if and only if $\lambda_1=\sqrt{\Delta}$ and
$\lambda_2=\cdots=\lambda_{n/2}$. Now the proof is complete.\qed

By expanding the right of Inequality (\ref{ELB1}), we obtain a
simplified lower bound, but it is a little weaker than the bound
(\ref{ELB1}).
\begin{cor}\label{ELBC}
Let $G^\sigma$ be a nonsingular oriented graph with order $n$,
maximum degree $\Delta$ and skew adjacency matrix $S$. Then
\begin{equation}\label{ELB2}
\mathcal{E}_s(G^\sigma)\geq 2\sqrt{\Delta}+n-2+\ln(\det(S))-\ln{\Delta}.
\end{equation}
Equality holds if and only if $G^\sigma$ is a union of $n/2$
disjoint arcs.
\end{cor}
\pf Note that $e^{x}\geq1+x$ for all $x$, where the equality holds
if and only if $x=0$. Combining the above inequality and Theorem
\ref{ELB}, we obtain that
\begin{eqnarray*}
\mathcal{E}_s(G^\sigma)&\geq& 2\sqrt{\Delta}+(n-2)\left(\frac{\det(S)}{\Delta}\right)^{\frac{1}{n-2}}\\
&=&2\sqrt{\Delta}+(n-2)e^\frac{\ln{(\det(S)/\Delta)}}{n-2}\\
&\geq&2\sqrt{\Delta}+(n-2)\left(1+\frac{\ln{(\det(S)/\Delta)}}{n-2}\right)\\
&=&2\sqrt{\Delta}+n-2+\ln(\det(S))-\ln{\Delta}\,\,.
\end{eqnarray*}
The equality holds in (\ref{ELB2}) if and only if all the
inequalities in the above consideration must be equalities, that is, 
$\lambda_1=\sqrt{\Delta}$, $\lambda_2=\cdots=\lambda_{n/2}$ and
$\det(S)=\Delta$.

It is easy to see that a union of $n/2$ disjoint arcs satisfies the
equality in (\ref{ELB2}). Conversely, if $\lambda_1=\sqrt{\Delta}$,
$\lambda_2=\cdots=\lambda_{n/2}$ and $\det(S)=\Delta$, then we find
that $\lambda_1=\sqrt{\Delta}$ and
$\lambda_2=\cdots=\lambda_{n/2}=1$. If $\Delta=1$, then $G^\sigma$
is a union of $n/2$ disjoint arcs. Suppose that $\Delta\geq 2$. Then
we get that $n\geq 4$. Considering the matrix $M=S^TS=(m_{ij})$,
where $S=(S_1,\ldots,S_n)$. By Theorem \ref{general}, we know that
$M$ is either a diagonal matrix with diagonal entries \{$\Delta,
\Delta, 1, \ldots, 1$\} or a matrix of form
$\left(\begin{array}{cc}\Delta & \mathbf{0}\\ \mathbf{0}^T &
M_1\end{array}\right)$.

If $M$ is a diagonal matrix, then the graph $G$ must contain two
vertices with maximum degree $\Delta$ and the other vertices with
degree 1. Since $\Delta\geq 2$, there must exist a path
$v_i\,v_j\,v_k$ satisfying $d_{G}(v_i)=1$. Then no matter how to
orient the graph $G$, we always get that $(S_i,S_k)\neq0$ and
$m_{ik}\neq 0$, which is a contradiction.

If $M$ is a matrix of form $\left(\begin{array}{cc}\Delta & \mathbf{0}\\
\mathbf{0}^T & M_1\end{array}\right)$, then $M_1$ must be symmetric
and have spectrum $\{\Delta,1,\ldots,1\}$. Denote by $A$ the vertex
set of the connected component of $G$ which contains the vertex
$v_1$. Suppose $A=\{v_1,v_2,\ldots,v_{t+1}\}$. Obviously, $t\geq
\Delta$.  By the spectral decomposition, we deduce that
$M_1=I_{n-1}+(\Delta-1)pp^T$, where $I_{n-1}$ is a unit matrix of
order $n-1$ and $p$ is a unit eigenvector of $M_1$ corresponding to
the eigenvalue $\Delta$. Therefore, we get
\begin{equation}\label{Form}
S^TS=\left(\begin{array}{cc}\Delta & \mathbf{0}\\ \mathbf{0}^T & I_{n-1}
+(\Delta-1)pp^T\end{array}\right).
\end{equation}
Suppose that $p=(p_1,p_2,\ldots,p_{n-1})^T$. Then we claim that $p$
satisfies following propositions:
\begin{enumerate}[(1)]
\item $p_1^2+\cdots+p_{n-1}^2=1$,
\item The degree sequence of $G$ is \{$\Delta, 1+(\Delta-1)p_1^2,\ldots,1+(\Delta-1)p_{n-1}^2$\},
\item $(S_{i+1},S_{j+1})=(\Delta-1)p_{i}p_{j}$ for any two distinct integer $1\leq i,j\leq n-1$,
\item $p_i\neq 0$ for $1\leq i\leq t$.
\end{enumerate}
The first proposition is trivial. The second and third propositions
follow from direct calculation of Equality (\ref{Form}). Now it
remains to prove the fourth proposition. If not, without loss of
generality, we suppose $p_1=0$. Then from the propositions (2) and
(3), we have $d_G(v_2)=1$ and $(S_2,S_j)=0$ for all $3\leq j\leq
t+1$, which is a contradiction.

From the above propositions, we observe that $(\Delta-1)p_{j}^2\geq
1$ for all $1\leq j\leq t$. Then
$\Delta-1=(\Delta-1)\sum_{j=1}^{n-1}p_j^2\geq t\geq \Delta$, which is a
contradiction.

Now we conclude that $G^\sigma$ must be a union of $n/2$ disjoint
arcs. The proof is thus complete.  \qed

\begin{figure}[h,t,b,p]
\begin{center}
\scalebox{1}[1]{\includegraphics{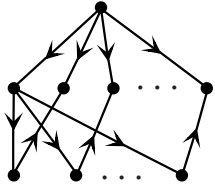}}
\end{center}
\caption{The oriented graph $G^\sigma$}\label{Fig2}
\end{figure}

For the general case, the bound (\ref{ELB2}) is not better than the
McClelland's lower bound (\ref{Mc}). For some cases, we find that
the bound (\ref{ELB2}) is better. For example, the oriented graph in
Figure \ref{Fig2} shows that the bound (\ref{ELB2}) is superior to
(\ref{Mc}). The oriented graph $G^\sigma$ has $n$ vertices, $3n/2-2$
arcs and maximum degree $n/2$. By calculation, we have
$\det(S)=n^2/4$.

Moreover, we find a class of oriented graphs which illustrate the
superiority of the bound (\ref{ELB2}). Let $\Gamma$ be the class of
connected oriented graphs of order $n\geq 600$, which satisfies the
following conditions:
\begin{equation}\label{class}
\frac{n}{2}\leq\Delta\leq \det(S)\leq n^2,\,\, m\leq 10n.
\end{equation}
Obviously, the oriented graph in Figure \ref{Fig2} belongs to $\Gamma$.

\begin{thm}\label{CLL}
The bound (\ref{ELB2}) is better than (\ref{Mc}) for any oriented
graph in $\Gamma$.
\end{thm}
\pf For any oriented graph $G^\sigma$ in $\Gamma$, we get that
\begin{equation*}
\ln{(\det(S))}-\ln{\Delta}=\ln\left(\frac{\det(S)}{\Delta}\right)\geq 0.
\end{equation*}
To prove the theorem, it is sufficient to prove that
\begin{equation}\label{equ}
2\sqrt{\Delta}+n-2\geq \sqrt{4m+n(n-2)(\det(S))^{\frac{2}{n}}}.
\end{equation}
Notice that $\ln(\det(S))\leq 2\ln{n}\leq \frac{n}{2}$ for $n\geq
600$. By Taylor's formula, we get
\begin{equation*}
(\det(S))^{\frac{2}{n}}=e^{\frac{2\ln(\det(S))}{n}}=1+\frac{2\ln(\det(S))}{n}+2e^{t_0}
\left(\frac{\ln(\det(S))}{n}\right)^2,
\end{equation*}
where the number $t_0$ satisfies that $0\leq t_0\leq
\frac{2\ln(\det(S))}{n}\leq \frac{4\ln{n}}{n}\leq1$. It follows that
$e^{t_0}\leq e$. We immediately obtain the following inequality:
\begin{equation}\label{taylor}
(\det(S))^{\frac{2}{n}}\leq 1+\frac{2\ln(\det(S))}{n}+2e\left(\frac{\ln(\det(S))}{n}\right)^2.
\end{equation}
To prove Inequality (\ref{equ}), we demonstrate that
\begin{eqnarray*}
&&4m+n(n-2)(\det(S))^{\frac{2}{n}}\\
&\leq& 4m+n(n-2)\left(1+\frac{2\ln(\det(S))}
{n}+2e\left(\frac{\ln(\det(S))}{n}\right)^2\right)\\
&\leq&38n+n^2+2(n-2)\ln(\det(S))+2e\left(\frac{n-2}{n}\right)
\left(\ln(\det(S))\right)^2\\
&\leq&38n+n^2+4(n-2)\ln(n)+8e\left(\frac{n-2}{n}\right)
\left(\ln(n)\right)^2\\
&\leq&n^2-2n+4+2\sqrt{2n}(n-2).
\end{eqnarray*}
The last inequality follows for $n\geq 600$. Note that
$(2\sqrt{\Delta}+n-2)^2\geq n^2-2n+4+2\sqrt{2n}(n-2)$. The proof is
thus complete. \qed

\noindent\textbf{Remark 4.2} By Theorem \ref{Mc1}, we know that the
bound (\ref{Mc}) is always superior to the McClelland's bound
obtained by Adiga et al. By Theorem \ref{ELB} and Corollary
\ref{ELBC}, the bound (\ref{ELB1}) is always superior to the bound
(\ref{ELB2}). For some cases, we obtain from Theorem \ref{CLL} that
the bound (\ref{ELB2}) is better than the bound (\ref{Mc}).

\end{document}